\documentclass[12pt]{amsart}
\usepackage{amsmath}
\usepackage{amscd}
\usepackage{graphics}
\usepackage{latexsym}

\oddsidemargin .in
\theoremstyle{plain}
\newtheorem{Thm}{Theorem}
\newtheorem{Lem}[Thm]{Lemma}

\newtheorem{Rm}{Remark}

\theoremstyle{definition}
\newtheorem{Def}{Definition}

\theoremstyle{remark}

\def\C{{\mathbb C}}
\def\A{{\mathbb A}}

\def\E{{\mathbb E}}
\def\R{{\mathbb R}}

\def\Z{{\mathbb Z}}
\def\O{{\mathbb O}}
\def\H{{\mathbb H}}
\def\V{{\mathbb V}}
\def\Q{{\mathbb Q}}

\def\CM{\mathcal M}
\def\CG{\mathcal G}

\def\CA{\mathcal A}

\def\CP{\mathcal P}
\def\CS{\mathcal S}
\def\CZ{\mathcal Z}

\def\s2x{\hbox{$S^2 \times S^2$}}

\def \Di {D\!\!\!\!/\,}

    \def\sqr#1#2{{\vcenter{\hrule height.#2pt
            \hbox{\vrule width.#2pt height#1pt \kern#1pt
            \vrule width.#2pt}\hrule height.#2pt}}}
    \def\square{\mathchoice\sqr67\sqr67\sqr{2.1}6\sqr{1.5}6}

\begin{document}

\title[]{ Associative submanifolds  of a $\mathbf{G_2}$ manifold\\   }
\author{Selman Akbulut and Sema Salur}
\thanks{First named author is partially supported by NSF grant DMS 9971440}
\address{Department  of Mathematics, Michigan State University,  MI, 48824}
\email{akbulut@math.msu.edu }
\address {Department. of Mathematics, Northwestern University, IL, 60208 }
\email{salur@math.northwestern.edu }
\date{\today}
\begin{abstract}
We  study deformations of  associative submanifolds $Y^3\subset M^7$  of a $G_2$ manifold $M^7$.  We show that the deformation space can be perturbed to  be smooth,   and it can be made compact and zero dimensional by constraining  it with  an additional equation. This allows us to associate local invariants to  associative submanifolds  of $M$. The local equations at each associative $Y$ are  restrictions of a global equation on a certain associated  Grassmann bundle over $ M$.

\end{abstract}
\maketitle

\setcounter{section}{-1}
\vspace{-0.4in}
\section{\bf Introduction}

 McLean  showed that, in a $G_2$ manifold $(M^7, \varphi)  $,   the space of associative submanifolds near a given one $Y^3,$  can be identified with  the  harmonic spinors on $Y$ twisted by a certain bundle $E$  (the kernel of a twisted Dirac operator)  \cite{m}. But  since we cannot control the cokernel of the Dirac operator (it has index zero), the  dimension of its  kernel might vary.  This is the obstruction to smoothness of the moduli space of associative submanifolds. This problem can be remedied either  by deforming the ambient $G_2$ structure (i.e. by deforming $\varphi $) or by deforming the connection in  the normal bundle \cite{as}. The first  process  might move $\varphi $ to a  non-integrable $G_2$ structure. If we are to view $ (M, \varphi)$  as an analogue of  a  symplectic manifold and $\varphi$ a symplectic form,  and view the associative submanifolds as  analogues of holomorphic curves, deforming  $\varphi$ would be too destructive. In the  second process we use the connections as auxiliary objects  to deform the associative submanifolds  in a larger space, just like deforming the holomorphic curves  by using almost complex structures (pseudo-holomorphic curves).  By  this approach we  obtain  the smoothness  of the moduli space. We get compactness  by relating the deformation equation to the Seiberg-Witten equations.

\vspace{.05in} 

 In this paper we summarize the results of  \cite{as} where we introduced {\it complex associative submanifolds} of $G_2$ manifolds;  they are associative submanifolds whose normal bundles carry  a $U(2)$ structure. This is no restriction, since every associative submanifold has this structure,  but if we require that their  deformations be compatible with the background connection we must have an integrability condition, i.e. the condition that  the connection  on their normal bundles  (induced by the $G_2$ background  metric)  reduce  to a $U(2)$ connection.  We call such manifolds {\it integrable complex associative submanifolds}. The prototype of these manifolds are holomorphic curves crossed by circles inside of $Calabi Yau  \times S^{1}$. We don't know how abundant these manifolds are in general.  One natural way to remove this restriction  is to use a generalized version of the Seiberg-Witten theory \cite{f},  the drawback of this is that it is a harder theory where we don't have automatic compactness. There is a  more general class of submanifolds, namely  the associative submanifolds whose normal bundles  carry a  $Spin^{c}(4)$ structure.  This structure allows us to break the normal bundle into pair of complex bundles, define necessary Dirac operator between them, and perform the deformations in these bundles via the standard Seiberg-Witten theory.  In this case we don't have a connection  compatibility problem,  because we perform the deformation in a larger space by using the whole  background connection with the  help of an auxiliary connection. More naturally,  we can define $U(2)$ or $Spin^{c}(4)$ structures on a $G_2$ manifold  itself, and induce them to all if its associative submanifolds.   We  can  also define global equations on the $15$-dimensional  associative Grassmannian bundle over a $G_2$ manifold $ (M,\varphi) $,  such a way that they restrict to the above mentioned Seiberg-Witten equations on  each of its associative submanifold $Y^3\subset M$.
 
  We would like to thank  G. Tian and R.Bryant for many  valuable comments,  and  also thank R. Kirby for constant encouragement.

\section{\bf Preliminaries}

Here we give some basic definitions and known facts concerning  the manifolds with special holonomy, the reader can find all these  interesting facts about them in \cite{b2},  \cite{b3}, \cite{hl}.
The exceptional Lie group $G_{2}$ can be defined as the subgroup of $GL(7,\R)$  fixing a particular $3$-form $\varphi_{0} \in \Omega^{3}(\R ^{7})$.
 Denote $e^{ijk}=dx^{i}\wedge dx^{j} \wedge dx^{k}\in \Omega^{3}(\R^7)$ , then
$$ G_{2}=\{ A \in GL(7,\R) \; | \; A^{*} \varphi_{0} =\varphi_{0}\; \} $$
\vspace{-0.2in}
$$\varphi _{0} =e^{123}+e^{145}+e^{167}+e^{246}-e^{257}-e^{347}-e^{356}$$
{\Def A smooth $7$-manifold $M^7$ has a {\it $G_{2}$ structure} if its tangent frame bundle reduces to a $G_{2}$ bundle. Equivalently
 $M^7$ has a {\it $G_{2}$ structure} if it has a 3-form $\varphi \in \Omega^{3}(M)$  such that  at each $x\in  M$  the pair $ (T_{x}(M), \varphi (x) )$ is  isomorphic to $(\R^{7}, \varphi_{0})$}.
 
 \vspace{.1in}

In fact $G_{2}$ is a $14$-dimensional subgroup of $SO(7)$. Since $GL(7,\R)$ acts on $\Omega^{3}(\R^7)$ with  stabilizer $G_{2}$, its orbit $\Omega^{3}_{+}(\R^7) $ is open for dimension reasons, so the choice of $\varphi_{0}$ in  the above definition is generic. In fact the action of  $GL^{+}(7,\R)$ has two orbits  containing $\pm \varphi_{0}$. The set of smooth $7$-manifolds with $G_{2}$-structures coincides with  the set of $7$-manifolds with spin structure, though this correspondence is not $1-1$. This is because $Spin(7)$ acts on $S^{7}$ with stabilizer $G_{2}$ inducing  the fibrations
$$G_{2}\to Spin(7)\to S^{7}\to BG_{2}\to BSpin(7)$$
and hence there is no obstruction to lifting maps  
$M^7\to BSpin(7)$ to $BG_{2}$.
The cotangent frame bundle $\CP^{*}(M)\to M$ of a manifold with $G_{2}$ structure $(M,\varphi)$  can be expressed as $\CP^{*}(M)=\cup_{x\in M} \;\CP^{*}_{x}(M)$, where
each fiber is given by:
\begin{equation}
\CP^{*}_{x}(M)=\{ u\in Hom (T_{x}(M), {\R}^{7})\;|\; u^{*}(\varphi_{0})=\varphi (x)\;\}
\end{equation}

 \vspace{.1in}

It turns out that any  $G_{2}$ structure $\varphi$ on $M^7$ gives an orientation $\mu \in \Omega^{7}(M)$ on $M$,   and in turn $\mu$  determines  a metric $g= \langle \;,\;\rangle$ on $M$, and a cross product structure $\times$  on its tangent bundle of $M$ as follows: Let  $i_{v}$ denote the interior product with a vector $v$, then 
\begin{equation}
g(u,v)=[ i_{u}(\varphi ) \wedge i_{v}(\varphi )\wedge \varphi  ]/\mu  
\end{equation}
\begin{equation}
\varphi (u,v,w)=g(u\times v,w)
\end{equation}

\noindent To emphasize the dependency on $\varphi $ sometimes $g$ is denoted by $g_{\varphi}$. There is a notion of a $G_2$ structure $\varphi $  on $M^7$ being integrable, which corresponds to $\varphi$ being harmonic.
{\Def A manifold with $G_{2}$ structure $(M,\varphi)$  is called a {\it $G_{2}$ manifold} if   the holonomy group of the Levi-Civita connection (of the metric $g_{\varphi}$) lies inside of  $G_2 $. Equivalently  $(M,\varphi)$ is a $G_{2}$ manifold if $d\varphi=d^{*}\varphi=0$ (i.e. $\varphi$ harmonic)}.

\vspace{.1in}

In short one can define a $G_{2}$ manifold to be any Riemannian manifold $(M^{7},g)$ whose holonomy group is contained in $G_{2}$;  then $\varphi $ and the cross product $\times$ come  as a consequence. It turns out that  the condition $\varphi$ being harmonic is equivalent to the condition that at each point $x\in M$ there is a chart  such that $\varphi (x)=\varphi_{0} + O(|x|^2)$ . 

For example if  $(X,\omega, \Omega)$ is a complex 3-dimensional Calabi-Yau manifold with 
K\"{a}hler form $\omega$ and a nowhere vanishing holomorphic 
3-form $\Omega$, then $X^6\times S^1$ has holonomy group 
$SU(3)\subset G_2$, hence  is a $G_2$ manifold. In this case
$\varphi$= Re $\Omega + \omega \wedge dt$.

\vspace{.05in}

{\Def  A 3-dimensional submanifold $Y$ of a $G_2$ manifold  $(M, \varphi )$ is
called {\em associative submanifold } if $\varphi|_Y\equiv vol(Y)$. This condition is equivalent to
$\chi|_Y\equiv 0$,  where $\chi \in \Omega^{3}(M, TM)$ is the  tangent bundle valued 3-form defined   by  the identity: }
\begin{equation}
\langle \chi (u,v,w) , z \rangle=*\varphi  (u,v,w,z)
\end{equation} 
The equivalence of these  conditions follows from  the `associator equality' of  \cite{hl}  $$\varphi  (u,v,w)^2 + |\chi (u,v,w)|^2/4= |u\wedge v\wedge w|^2$$ We will denote the tangent bundle valued $3$-form in $\R^7$ corresponding to $\varphi_{0}$ by $\chi_{0}$.
Throughout this paper we will denote the sections of a bundle $\xi \to Y$ by $\Omega^{0}(Y,\xi)$ or simply by  $\Omega^{0}(\xi)$, and the bundle valued $p$-forms by $\Omega^{p}(\xi)=\Omega^{0}(\Lambda^{p}T^{*}Y\otimes \xi)$.  
We will denote the  coframe bundle by  $\CP^{*}(M)\to M$ and its adapted frame bundle by $\CP(M)$. They can be $G_{2}$ or $SO(7)$ frame bundles;  when needed we will specify them by the notations  $\CP_{SO(7)}(M)$ or  $\CP_{G_{2}}(M)$.

\section{ \bf  Grassmann bundles}

Let $G(3,7)$ be the Grassmann manifold  of oriented $3$-planes in $\R^7$. Let $M^7$ be any smooth $7$-manifold, and let  $\tilde{M}$ be the $3$-planes in $T(M)$, i.e.  $\tilde{M}\to M$ is the bundle of Grassmannians over $M$ defined by
$$ \tilde{M}= \CP_{SO(7)}(M)\times_{SO(7)} G(3,7) \to M$$

\noindent Let $\xi \to G(3,7)$ be the universal $\R^3$ bundle, and $\nu =\xi^{\perp}\to G(3,7)$ be the dual $\R^{4}$ bundle. Hence
$Hom(\xi, \nu)=\xi^{*}\otimes \nu \longrightarrow G(3,7)$ is just  the tangent bundle  $TG(3,7)$.
$\xi$, $\nu$ extend fiberwise to give  bundles $\Xi \to  \tilde{M}$, $\V \to \tilde{M}$ respectively, and let $\Xi^*$ denote the dual of $\Xi$. Notice   that 
$ Hom (\Xi, \V)= \Xi^{*}\otimes \V \to \tilde{M}\; $
is just the bundle  of vertical vectors  $\; T_{vert}(\tilde{M}) $ of $T(\tilde{M}) \to M$, i.e. it is the bundle of tangent vectors to the fibers of  $\pi: \tilde{M}\to M$. Hence 
\begin{equation} 
T\tilde{M}=  T_{vert}(\tilde{M}) \oplus  \Xi \oplus \V  = ( \Xi^{*}\otimes \V ) \oplus  \pi^{*} TM
\end{equation}  

Let  $\CP(\V)\to \tilde{M}$ be the $SO(4)$ frame bundle of the vector bundle $\V$, and identify $\R^4$ with  the quaternions $ \H $, and identify $SU(2)$ with the unit quaternions $Sp(1)$. Recall  that $SO(4)$ is the equivalence classes of pairs $[\;q,\lambda \;]$ of unit quaternions
$$SO(4)=(SU(2)\times SU(2))/\Z_{2}$$ Hence   $\V\to \tilde{M} $  is the associated vector bundle to $\CP(\V)$ via  the $SO(4)$ representation 
\begin{equation}
x \mapsto  qx\lambda^{-1}
\end{equation}
 There is  a pair of $\R^3=im(\H)$ bundles over $\tilde{M}$ corresponding to the left and right $SO(3)$ reductions of $SO(4)$,  so  they are given by the $SO(3)$ representations
\begin{equation}
\begin{array}{lcc}
 \lambda_{+}(\V)\;:  &x  \mapsto qx\ q^{-1} &   \\
\lambda_{-}(\V)  \;: &y \mapsto  \lambda y  \lambda ^{-1} & 
 \end{array}
 \end{equation}
The map $x\otimes y\mapsto xy$  gives  actions $\lambda_{+}(\V)\otimes \V \to \V  \;\; \hbox {and}\;\; \V\otimes \lambda_{-}(\V) \to \V$;   by combining we can think of it as one conjugation  action 
\begin{equation}
(\lambda_{+}(\V)\otimes \lambda_{-}(\V))\otimes  \V \to \V 
\end{equation}

\vspace{.1in}

If  the $SO(4)$ bundle $\V \to \tilde{M}$ lifts to a $Spin(4)= SU(2)\times SU(2)$ bundle (locally it does),  we get two additional bundles  over $\tilde{M}$
 \begin{equation}
\begin{array}{lcc}
\;\CS\;: & y \mapsto  q y  & \\
\;\E\;: & \; \; y \mapsto  y  \lambda ^{-1} & \\
 \end{array}
 \end{equation}
 
\noindent which gives $\V$  as a  tensor product of two quaternionic line bundles  $\V= \CS \otimes_{\H} \E$.  In particular $\lambda_{+}(\V)=ad(\CS)$ and $\lambda_{-}(\V)= ad(\E)$, i.e. they are the $SO(3)$ reductions of the $SU(2)$ bundles $\CS$ and  $\E$.  In particular there is  multiplication map 
$ \CS \otimes \E \to \V$.

\section{\bf  Associative Grassmann bundles}

  Now consider  the {\it Grassmannian of associative $3$-planes} $G^{\varphi}(3,7) $  in $\R^7 $, consisting of elements $L\in  G(3,7)$ with the property $\varphi_{0}|_{L}=vol(L)$ (or equivalently $\chi_{0}|_{L}=0$). $G_2$ acts on $G^{\varphi}(3,7) $  transitively  with  the stabilizer $SO(4)$, giving the the identification $G^{\varphi}(3,7)=G_{2}/SO(4) $.  Identify  the imaginary octonions $\R^7 =\mbox{Im} (\O )\cong im (\H)\oplus \H $,  then the action of the  subgroup $SO(4)\subset G_{2}$ on $\R^7$ is given by:  
\begin{equation}
\left(
\begin{array}{cc}
  A   &  0 \\
   0  &  \rho(A) \\  
\end{array}
\right)
\end{equation}
where $\rho :  SO(4) = (SU(2)\times SU(2))/ \Z_{2} \to SO(3)$ is the projection of the first factor (\cite{hl}).  Now let  $M^7$  be a $G_2$ manifold.  As above,  this time we can construct  the bundle of associative Grassmannians over $M$
$$\tilde{M}_{\varphi} = \CP_{G_{2}}(M) \times _{G_{2}}G^{\varphi}(3,7) \to M$$ 
which is just the quotient bundle $\tilde {M}_{\varphi}=\CP(M)/SO(4)\longrightarrow \CP(M)/G_{2} =M$. As in the previous section,   the universal bundles $\xi, \;\nu=\xi^{\perp} \to G^{\varphi}(3,7)$ induce $3$  and $4$ plane bundles $\Xi\to\tilde{ M}_{\varphi}$ and $\V\to \tilde{M}_{\varphi}$ (i.e. restricting the universal  bundles  from $ \tilde{M}$).   Also 

\begin{equation} 
T\tilde{M}_{\varphi }=T_{vert } (\tilde{M}_{\varphi }) \oplus  \Xi \oplus \V  
\end{equation}

So,  in  this  associative case,  we have an important identification $\Xi= \lambda_{+}(\V)$  (as bundles over $\tilde{M}_{\varphi}$). For simplicity we will denote   $\Q= \lambda_{-}(\V)$.  As before we have an  action $ \V \otimes \Q \to \V$, and the similar action of $\lambda_{+}(\V)$ on $\V$.  We will break the dual of this last  action into a pair of 
Clifford multiplications on $\V^{\pm}=\V$
 \begin{equation}
 \Xi^{*} \otimes \V^{\pm} \to \V^{\mp}  
\end{equation}
 given by   $x\otimes y \mapsto -\bar{x}y$ and $x\otimes y \mapsto xy$,  on $\V^{+}$ and $\V^{-}$ respectively. By writing $x_{1}\wedge x_{2} =\frac{1}{2}(x_{1}\otimes x_{2}-x_{2}\otimes x_{1})$  we can extend this  to an action  
of the whole exterior algebra $\Lambda^{*} (\Xi^{*})$. For example  $ ( x_{1} \wedge x_{2}) \otimes y \mapsto \frac{1}{2}  (-x_{1} \bar{x}_{2}+x_{2}\bar{x}_{1} )  y =Im(x_{2}\bar{x}_{1}) y $ gives 
\begin{equation}
 \Lambda^{2}(\Xi^{*})\otimes \V  \to \V
 \end{equation}
    
  \vspace{.1in}

 Recall $T_{vert}(\tilde{M})=  \Xi^{*}\otimes \V $ is the  the subbundle of  vertical  vectors of  $T (\tilde {M}) \to M$.   The total space $E(\nu_{\varphi})$ of the normal bundle of the imbedding $ \tilde {M}_{\varphi } \subset \tilde{M}$  should be thought of an open tubular neighborhood of $ \tilde {M}_{\varphi }$ in $\tilde{M}$, and it has a nice description:  

 \vspace{.05in}

{\Lem  (\cite{m})  The normal bundle $\nu_{\varphi}$ of $ \tilde {M}_{\varphi } \subset \tilde{M}$  is isomorphic to $\V$,  and the bundle of vertical vectors $T_{vert} (\tilde{M}_{\varphi })$  of $T(\tilde{M}_{\varphi }) \to M$  is the kernel of the Clifford multiplication $ c: \Xi ^{*} \otimes \V \to \V$. }

 \vspace{.1in}
 
  Hence  we have $ \; \Xi^{*} \otimes \V| _{{\tilde{M}_{\varphi}}}  =T_{vert }(\tilde{M}_{\varphi })\oplus \nu _{\varphi }$;  from the  exact sequence over $\tilde{M}_{\varphi }$
\begin{equation}
T_{vert} (\tilde{M}_{\varphi })\to \Xi^{*}\otimes \V|_{{\tilde{M}_{\varphi}}} \stackrel{c}{\longrightarrow} \V|_{\tilde{M}_{\varphi }}
\end{equation}

\noindent the  quotient bundle, $ T _{vert} (\tilde{M}) / T_{vert} (\tilde{M}_{\varphi })$  is isomorphic to $\nu_{\varphi}$.  Also, 
 if $ \{e^{i}\}$ is a local orthonormal  basis for $\Xi^{*} $,  then  the projection $ \pi_{\varphi} :T_{vert} (\tilde{M})|_{\tilde{M}_{\varphi}} \to T_{vert} (\tilde{M}_{\varphi}) $ is  given by 
\begin{equation}
\pi_{\varphi}(a\otimes v)=a\otimes v +(1/3)\sum _{j=1}^{3} e^{j}\otimes e^{j}(a.v)
\end{equation}

{\Rm The lemma holds since  the Lie algebra inclusion $ g_2\subset so(7) $ is given by
\[
\left(
\begin{array}{ccc}
a  &  \beta    \\
-\beta^t  &\rho ( a )    \\
 \end{array}
\right)
\]
\noindent  where $a\in so(4)$  is $y\mapsto qy-y\lambda$, and  $\rho(a)\in so(3)$ is  
 $x\mapsto qx-xq$. So the tangent space inclusion of  $G_{2}/SO(4)\subset SO(7)/SO(4)\times SO(3)$ is given by the matrix $\beta$. And if we write $\beta$ as column vectors $\beta=(\beta_1, \beta_2, \beta_3)$, then $\beta_{1}i+\beta_{2}j+\beta_{3}k=0$ (\cite{m},  \cite{mc})}

\section {\bf Associative submanifolds}

 Any imbedding of a $3$-manifold   $f: Y^{3}\hookrightarrow M^7 $  induces a bundle $\tilde{Y}\to Y$ with fibers $G(3,7)$, and the Gauss map of $f$ canonically lifts to  an imbedding  $\tilde{f}: Y \hookrightarrow   \tilde{M}$: 
\begin{equation}
\begin{array}{lcl}
  & & \tilde{M}\supset \tilde{M}_{\varphi } \\
  \hspace{.25in} \tilde{f}   \hspace{-.1in}&  \nearrow \;  &   \downarrow\\
\;Y\;\;& \stackrel{f}{\longrightarrow}   & M
\end{array}
\end{equation}
 Also,  the pull-backs $\tilde{f}^{*} \Xi=T (Y)$ and  $\tilde{f}^{*} \V=\nu(Y)$ give the tangent and normal bundles of $Y$.  Furthermore,   if $f$ is an imbedding as an  associative submanifold in  a $G_2$ manifold  $(M, \varphi)$,  then the image of $\tilde{f}$ lands in $ \tilde {M}_{\varphi }$.  We will denote this canonical  lifting of any $3$-manifold $Y\subset M$ by $\tilde{Y}\subset \tilde{M}$.

\vspace{.05in}
$ \tilde {M}_{\varphi }$ can be thought of as a universal space parametrizing associative submanifolds of $M$. In particular,  if $\tilde{f}: Y\hookrightarrow  M_{\varphi}$ is the lifting of an associative submanifold, by pulling back we see that  the principal  $SO(4)$   bundle $\CP(\V)\to \tilde{M}_{\varphi} $ induces an $SO(4)$-bundle $\CP (Y)\to Y$, and gives  the  following  vector bundles over $Y$ via the  representations:  
\begin{equation}
\begin{array}{lcc}
\; \nu(Y) \hspace{.2in}:   &y \mapsto  q y\lambda^{-1}   & \hspace{1in}  \\
 \; T (Y)  \hspace{.1in}:  & x \mapsto  qx\ q^{-1} &   \\
\;Q  \hspace{.4in}: & y  \mapsto  \lambda y  \lambda ^{-1} & \\
 \end{array}
 \end{equation}
where $[q,\lambda]\in SO(4)$. Note that $T(Y)= \lambda_{+}(\nu ) $ and $Q= \lambda_{-}(\nu )$, where $\nu=\nu(Y)$. Also we can identify $T^{*}Y$ with $TY$ by the induced metric. We  have the actions 
$T^{*}Y \otimes \nu \to \nu  $,    
$\Lambda^{*}(T^{*}Y) \otimes \nu \to \nu  $,  $ \nu \otimes Q \to \nu$ from above. 
 By combining  them we can define  an action of  the bundle valued  differential forms 
\begin{equation}
\Omega^{*}(Y;Q) \otimes\Omega^{0}(Y,\nu) \to \Omega^{0}(Y,\nu)
\end{equation}

\noindent For example, if $F=  (x_{1} \wedge x_{2})\otimes y$, and $y\in \Omega^{0}(Y,Q)$, and $z  \in \Omega^{0}(Y, \nu)$ then
\begin{equation}
F\otimes z \mapsto Im (x_{2}\bar{x}_{1}) z y 
\end{equation}

\vspace{.1in}

Let $\L=\Lambda^{3}(\Xi)\to \tilde{M}$  be the determinant (real) line  bundle  . Note that the  definition (4) implies that  $\chi $ maps every oriented $3$-plane in $T_{x}(M)$ to its complementary subspace,  so $\chi$ gives a bundle map $\L \to \V$ over $\tilde{M}$, which is a section of $\L^{*}\otimes \V \to \tilde{M}.$  Also, $\Xi$ being an oriented bundle implies $\L$ is trivial  so $ \chi $ actually gives  a section
\begin{equation}
 \chi  \in \Omega^{0}(\tilde{M},  \V)
 \end{equation}
Clearly $ \tilde {M}_{\varphi }\subset \tilde{M}$  is the codimension $4$  submanifold given  as the zeros of this section. So associative submanifolds $Y\subset M$ are characterized by the condition $\chi|_{\tilde{Y}}=0$,  where $\tilde{Y}\subset \tilde{M}$ is the canonical lifting of $Y$. Similarly $\varphi $ defines a map $\varphi: \tilde{M}\to \R$.

 \section{\bf Dirac operator }

 In general,  the normal bundle $\nu=\nu(Y)$  of any orientable $3$-manifold $Y$  in a $G_2$ manifold $ (M,\varphi) $ has a Spin(4) structure  (e.g. \cite{b2}). Hence we have $SU(2)$ bundles $S$ and $E$   over  $Y$ (the pull-backs of $\CS$ and $\E $ in (9)) such that $\nu =S\otimes_{\H} E$,  with $SO(3)$ reductions $ad \;S=\lambda_{+} (\nu)$ and  $ad \;E= \lambda_{-}(\nu)$. Note that  $ad \;E$  is also  the bundle of endomorphisms $End (E)$.  If  $Y$ is  associative,  then  the $SU(2)$ bundle $S$  reduces  to $TY,$  i.e. $S$ becomes  the spinor bundle of $Y$,  and so $\nu $ becomes  a twisted spinor bundle on $Y$. 

\vspace{.1in} 

The Levi-Civita connection of the $G_{2}$ metric of $(M,\varphi) $ induces  connections on  the associated bundles $\V$ and $\Xi$ on $\tilde{M}$. In particular it induces connections on the tangent and normal bundles of any submanifold $Y^{3}\subset M$. We will call  these connections  the {\it background connections}. Let $\A_{0}$ be the induced connection on  the normal bundle $\nu=S\otimes E$. From the   Lie algebra decomposition  $ so(4)=so(3)\oplus so(3) $,  we can write $ \A_{0}= A_{0}\oplus B_{0}$, where $A_{0}$ and  $B_{0}$ are connections on $S$ and $E$ respectively. 

\vspace{.05in}

 Let   $\CA(S)$  and $\CA(E)$ be the  set of connections on  the  bundles $S$  and $E$. Hence   connections $ A\in \CA(S) $,  $ B\in \CA(E)$  are in the form $A=A_{0} +a$, $B=B_{0}+b$,   where $a\in  \Omega^{1}(Y, ad \;S) $ and $b\in  \Omega^{1}(Y, ad\; E) $.  So    $\Omega^{1}(Y, \lambda_{\pm}(\nu))$ parametrizes connections on $S$ and $E$, and the connections on $\nu $ are in the form $\A= A\oplus B$. To emphasize the dependency on $a$ and $b$ we will sometimes denote $\A=\A(a,b)$, and $\A_{0}=\A(0,0)$.

 \vspace{.1in}
 
Now, let $Y^3\subset M$ be any smooth  manifold.  We can express  the covariant derivative $\nabla_{\A}: \Omega^{0}(Y, \nu)\to \Omega^{1}(Y, \nu)$ on $\nu$  by $ \nabla_{A}=\sum e^i \;\otimes \nabla_{e_i}$, where $\{e_{i}\}$ and $\{e^{i}\}$ are  orthonormal tangent and cotangent  frame fields of $Y$,  respectively. Furthermore if $Y$ is an associative submanifold,  we can use the Clifford multiplication  (12) to form the twisted Dirac operator  $\Di _{A}: \Omega^{0}(Y,\nu) \to \Omega^{0}(Y,\nu)$ 
  \begin{equation}
  \Di _{\A}=\sum e^i \;. \nabla_{e_i}
 \end{equation}
The sections lying in  the kernel of this operator are usually called  harmonic spinors twisted by $(E, \A)$. Elements of  the kernel of $\Di_{A_{0}}$ are called the  harmonic spinors twisted by $E$,  or just the twisted harmonic spinors.

\section {\bf Complex associative submanifolds}

There is an interesting class of associative submanifolds which we  call {\it complex associative}, they are defined as follows: The  subgroups $U(2)\subset SO(4)\subset G_{2}$,  from  
$$  (SU(2)\times S^1)/ \Z_{2}\subset (SU(2)\times SU(2))/\Z_{2}\subset G_{2}$$
  
\noindent  give a $U(2)$-principal  bundle  
 $\CP_{G_2} (M)\to  \bar{M}_{\varphi }=\CP_{G_2}(M) /U(2) $. Note that $ \bar{M}_{\varphi }$ is   the total space of  an $S^{2}$ bundle  $ \bar{M}_{\varphi}\to \tilde{M}_{\varphi }=\CP(M)/SO(4)$.  This is just  the  sphere bundle of the $\R^{3}$-bundle $\lambda_{-}(\V)\to \tilde{M}_{\varphi}$. So we can identify the unit sections  $j\in \lambda_{-}(\V)$ with complex structures on $\V$, i.e. the right reductions of $SO(4)$ to $U(2)$, where
$$U(2)=\mbox{Spin}^{c}(3)=(SU(2) \times S^1) /\Z_{2}$$

Just as $SO(4)$ is the  stabilizer subgroup of the action of $G_{2}$ on $G^{\varphi }(3,7)$, the subgroup $U(2)\subset SO(4)$ is the stabilizer of the action of $G_{2}$ on the "framed" version of $G^{\varphi }(3,7)$. That is, if 
$\nu\to G^{\varphi}(3,7) $ is the 
dual  universal bundle, then the  action of $G_{2}$ on $ G^{\varphi}(3,7)$ extends to an action to the sphere bundle of  $\lambda_{-}(\nu)$ with stabilizer $U(2)$.
\vspace{.1in}

By the  representations  of $U(2)$ induced from $SO(4)$ we can form  the same vector bundles $\V, \Xi $ on $\bar{M}_{\varphi }$ as in  Sections (2)  and (3).  
Except in this case we have $\lambda\in S^{1}$ and hence $\lambda_{-}(\V)$ becomes the trivial bundle,  and  $\V$ becomes a $\C^2$ bundle.   The reduction $j$ gives an extra line bundle $\L$ via the representation  $\rho_{2}(q,\lambda)(y)= y \lambda^{-2} $
(the determinant line bundle of $\V$). Similarly we have the same actions of (12), (13).  Since $\lambda_{-}(\V)$ is trivial its action becomes the multiplication with complex numbers.
In this case we have a useful quadratic bundle map $\sigma : \V \otimes \V \to \Lambda^{2}(\Xi^{*}) _{\C}$ given by  
\begin{equation}
 \sigma(x,y)=-\frac{1}{2}(x i \bar{y})i
\end{equation}
In particular if $x=z+jw \in \O$ (octonions),  we have $\sigma(x,x)= \frac{|z|^2 -|w|^2}{2}+ j \bar{z} w $.
  \vspace{.1in}
  
So if   $f:Y^{3}\hookrightarrow M$ is an associative submanifold, and $\tilde{f} : Y^{3}\hookrightarrow \tilde{M}_{\varphi }$ is the lifting of its Gauss map, then  the normal bundle of $f$ has a $U(2)$ structure  if and only if $\tilde{f}$ lifts   
\begin{equation}
\begin{array}{rrc}
  & & \bar{M}_{\varphi}  \\
 \bar{f} \hspace{-.1in} &  \nearrow \;  &   \downarrow\\
\;\; Y\;\; & \stackrel{\tilde{f}}{\longrightarrow}   &\tilde{ M}_{\varphi}
\end{array}
\end{equation}

As before we get $\bar{f}^{*} \Xi=T(Y)$ and  $\bar{f}^{*} \V=\nu(Y)$,  plus a line bundle $L=\bar{f}^{*}(\L)$, and $\nu$ becomes a complex $U(2)$ bundle, which we will call $W$. So $W$ is  the spinor bundle $W(L)\to Y$ of  a $Spin^c (3) $ structure on $Y$, and $L$ is its determinant bundle.  If  this $Spin^{c}(3)$ structure comes from a $Spin(3)$ structure,  then $W=S\otimes L^{1/2}$ where $S$ is the spinor  bundle on $Y$. So in concrete terms:

{\Def  A {\it complex associative submanifold} $(Y,j )$ of a $G_2$ manifold $M$ is an associative manifold in $M$ with a unit section $j\in \Omega^{0}(\lambda_{-}(\nu ))$}, where $\nu=\nu(Y)$.

\vspace{.1in}
 Note that the normal bundle $\nu(Y)$ of any associative submanifold $Y$ has a $U(2)$ structure (in fact it is trivial \cite{bsh}), so  liftings to $\tilde{M}_{\varphi}$ always exists. We will call any $Y$ with such a lifting a  {\it complex associative submanifold}.  But for such a lifting the background $SO(4)$ connection $\A_{0}$ on $\nu=\nu(Y)$ may not reduce to  a $U(2)$ connection.  Let  $j \in \lambda_{-}(\nu )$ be the unit section describing  the $U(2)$ structure on $\nu $. Then  $\A_{0}$ reduces to a $U(2)$ connection if $\nabla_{A_{0}}(j)=0$, where $A_{0}$ is the connection on  $ \lambda_{-}(\nu )$ induced from
 the background connection $\A_{0}$, and $\nabla_{A_{0}}$ is the covariant differentiation.
\begin{equation}
 \Omega^{0}(Y,  \lambda_{-}(\nu )) \stackrel {\nabla_{A_{0}}}
 {\longrightarrow }
  \Omega^{1}( Y, \lambda_{-}(\nu )) 
 \end{equation}
 
{\Def  An  {\it integrable complex associative submanifold} $(Y,j )$ of a $G_2$ manifold $M$ is an  associative manifold in $M$ with a unit section $j\in \Omega^{0}(\lambda_{-}(\nu ))$ with  $\nabla_{A_{0}}(j)=0$.} 

\vspace{.1in}

 A natural way to obtain complex associative submanifolds is by pulling back a universal complex structure from the ambient $G_2$ manifold $(M,\varphi)$:
 
 {\Def $(M,\varphi, J)$ is called a {\it complex $G_2$ manifold} if $J:\tilde{M}_{\varphi } \to \bar{M}_{\varphi} $. is a lifting}

\vspace{.05in}

 If $Y$ be an integrable complex associative submanifold of a $G_2$ manifold $M$, and $\A_{0}$ be the background $U(2)$ connection on its normal bundle $W$ induced by the $G_2$ metric. Then  $\A_{0}=S_{0}\oplus A_{0}$ where $S_{0}$ is a connection on $\lambda_{+}(W)=T(Y)$  and $A_{0}$ is a connection on the line bundle $L$. We can deform the background connection $\A_{0}=\A(0)\mapsto \A(a) $,  where $\A(a)=S_{0} \oplus  (A_{0} +a) $, where $a\in \Omega^{1}(Y)$  (e.g. Sec 8). 
  
 {\Rm Recall that if $(X,\omega, \Omega)$ is a complex 3-dimensional Calabi-Yau manifold  then $(X^6\times S^1, \varphi )$ is a $G_2$ manifold, where
$\varphi$= Re $\Omega + \omega \wedge dt$. The submanifolds of the form $C\times S^1$ where $C\subset X$ a holomorphic curve,  are examples of integrable complex associative submanifolds of $X\times S^1$. }

 \section{\bf Deformations }

  In \cite{m}   McLean  showed that,  in a neighborhood of an associative submanifold $Y$ in a  $G_2$ manifold $(M,\varphi) $,  the space of associative submanifolds can be identified with  the  harmonic spinors on $Y$ twisted by $E$. Since the  cokernel of the Dirac operator can vary, we cannot determine the dimension of its kernel  (it always  has zero index on $Y^3$). We will remedy this by deforming $Y$  in a larger class of submanifolds.
 To motivate our approach  let us sketch a proof of  McLean's theorem by adapting  the  explanation in \cite{b3}. 
Let  $Y\subset M$ be an associative submanifold, $Y$ will determine a lifting $\tilde{Y} \subset \tilde{M}_{\varphi }$. Let us recall that the  $G_2$ structure $\varphi $ gives a metric connection on $M$, hence it gives a covariant differentiation  in  the normal bundle $\nu(Y)=\nu$
$$\nabla_{\A_{0}}: \Omega^{0}(Y, \nu)\to \Omega^{1}(Y, \nu)=\Omega^{0}(Y, \;T^{*}Y\otimes \nu)$$
Recall that we  identified  $T^{*}_{y }(Y) \otimes \nu_{y}(Y)$ with  a copy of $TG(3,7)$, more precisely we identified it with  the tangent space of the Grasmannians of $3$-planes in $T_{y}(M)$. So the  covariant derivative lifts normal vector fields $v$  of $Y\subset M$  to vertical vector fields $\tilde{v}$  in $T(\tilde{M})|_{\tilde{Y}}$. We want the normal vector fields $v$ of $Y$ to always move $Y$ in the class of associative submanifolds of  $M$,  i.e. we want the liftings $\tilde{Y}_{v}$ of the nearby copies $Y_{v}$  of $Y$  (pushed off by the vector field $v$) to  lie in  
$T(\tilde{M}_{\varphi })\subset T(\tilde{M}) $ upstairs, i.e. we want the component of $\tilde{v}$  in the direction  of  the normal bundle $\tilde{M}_{\varphi}\subset \tilde{M}$  to vanish. By Lemma 1,  this means $\nabla_{\A_{0}}(v)$ should  be in the kernel of the Clifford multiplication $ \; c: \Omega^{0}(T^{*}(Y)\otimes \nu)\to \Omega^{0}(\nu)$, i.e. 
\begin{equation}
\Di_{\A_{0}}(v)=c (\nabla_{\A_{0}}(v))=0
\end{equation}

  \begin{figure}[ht]  \begin{center}  
\includegraphics{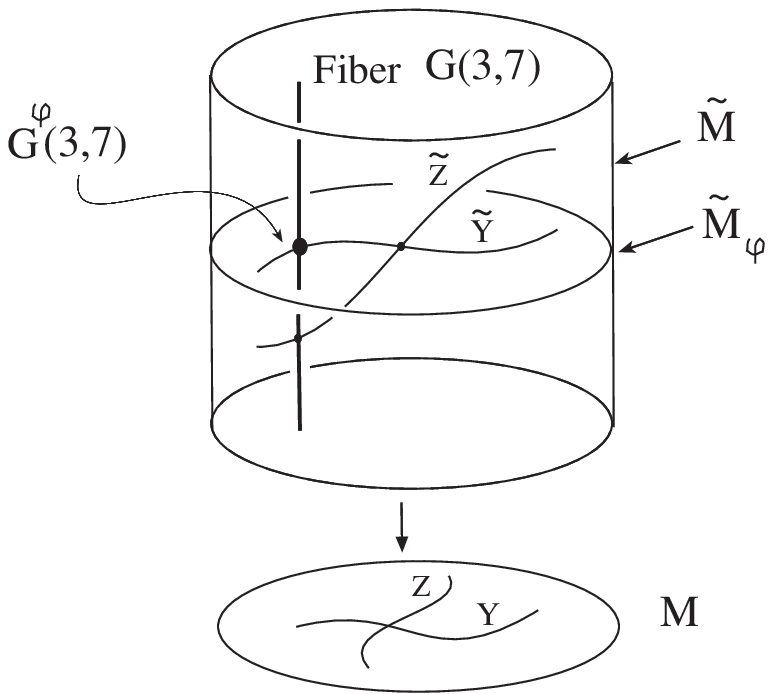}   \caption{}    \end{center}
   \end{figure}

 \noindent Here $\Di_{\A_{0}}$ is the Dirac operator induced by the  background connection $\A_{0}$, which is the composition 
\begin{equation}
\Omega^{0}(Y, \nu) \stackrel{\nabla_{\A_{0}}}{\longrightarrow } \Omega^{0}(Y, T^{*}Y\otimes 
\nu)\stackrel{c}{\to}  \Omega^{0}(Y, \nu) 
\end{equation}

To sum up what we have done: Let   $f: Y^{3}\hookrightarrow M$ be  an associative manifold in a  $G_2$ manifold, and $\tilde{f}: Y^{3}\hookrightarrow \tilde{M}_{\varphi }$ be its canonical lifting. By using the metric connection, the map on the space of immersions  $\Phi: Im (Y,M) \to Im (Y,\tilde{M}) $ given by  $f\mapsto \tilde{f}$  induces  a map on the tangent spaces 
\begin{equation}
d\Phi : T_{f} Im (Y,M) \longrightarrow T_{\tilde{f}} Im (Y,\tilde{M})\supset  T_{\tilde{f}} Im (Y,\tilde{M}_{\varphi  })
\end{equation}
and the condition that $d\Phi (v)$ lie in  $T_{\tilde{f}} Im (Y,\tilde{M}_{\varphi  })$ is given by  $\Di_{\A_{0}}(v)=0$. 

\vspace{.05in}

Now  $\Phi $ may not be transversal to $Im(Y,\tilde{M}_{\varphi })$,  i.e.  the cokernel of $\Di_{\A_{0}}$ may be non-zero, but we can make it zero by deforming  the connection in the normal bundle (i.e. in  $E$ or $S$  as indicated in Section 5)  $ \A_{0}=\A(0,0) \to A= \A(a,b) $, where $a\in \Omega^{1}(Y,\lambda_{+}(\nu ))$ and  $b\in \Omega^{1}(Y,\lambda_{-}(\nu ))$, and $\nu =\nu(Y) $. Then  (26) becomes parametrized  
\begin{equation}
\Omega^{0}(\nu) \times  \Omega^{1}(\lambda_{\pm}(\nu ))  \stackrel{\Di_{\A}}{\longrightarrow} \Omega^{0}(\nu) 
 \end{equation}
that is,   (25)  becomes a  twisted Dirac equation  $\Di_{\A}(v)= c(\nabla _{\A}(v))=\Di_{\A_{0}}(v) + \alpha v=0 $, where $\alpha=(a,b)$. This gives a  smooth  moduli space parametrized by $\Omega^{0}(\lambda _{\pm}(\nu))$. Then by a generic choice of $\alpha \in \Omega^{0}(\lambda _{\pm}(\nu))$ we get a zero-dimensional perturbed moduli space, whose elements perhaps should be called $A$ -associative submanifolds  (to make an  analogy with $J$-holomorphic curves). 

 {\Thm In a $G_2$ manifold $(M,\varphi )$,  for a generic  $A$, the space of $A$-associative submanifolds
in a neigborhood of an associative manifold  is a  zero dimensional oriented smooth manifold}.

\vspace{.1in}

In \cite{as} it was also shown that by deforming the $G_2$ structure $\varphi$ we can also obtain smoothness. To get compactness we will couple the  Dirac equation with another equation to get Seiberg-Witten  equations. We will first treat the complex associatives.

 \section{\bf Deforming Complex Associative Submanifolds }

Let  $\CM_{\C}(M, \varphi ) $ be the set of  complex associative submanifolds $(Y,j)$  of a $G_2$ manifold $(M,\varphi )$.  We wish  to study the local structure of  $\CM_{\C }(M,\varphi ) $ near a particular   $Y$.  The normal bundle of $Y$ is a $U(2)$  bundle  $\nu(Y)= W$,  and $L\to Y$ is its determinant  line bundle.  If $(Y, j) $ is integrable we may assume $\A_{0} =S_{0}\oplus A_{0 } $ is  the background  connection on $W$, where $S_{0}$ is a connection on $\lambda_{+}(W)=T(Y)$,  and $A_{0}$ is a connection on $L$.  As before we will denote  $\A= \A (a) = S_{0} \oplus  (A_{0} +a) $,  where  $a\in \Omega^{1}(Y) =T_{\A_{0}}\CA(L)$ (the tangent space to the space of connections).  Now (28) becomes
\begin{equation}
  \Omega^{0}(Y, W) \times \Omega^{1}(Y,i\R)  \stackrel{\Di_{\A}}{\longrightarrow} \Omega^{0}(Y, W) 
 \end{equation}
which can be thought of  the derivative of a  similarly defined map 
\begin{equation}
     \Omega^{0}(Y, W) \times \CA(\L)   \to  \Omega^{0}(Y, W)
  \end{equation}
  
  \begin{figure}[ht]  \begin{center}  
\includegraphics{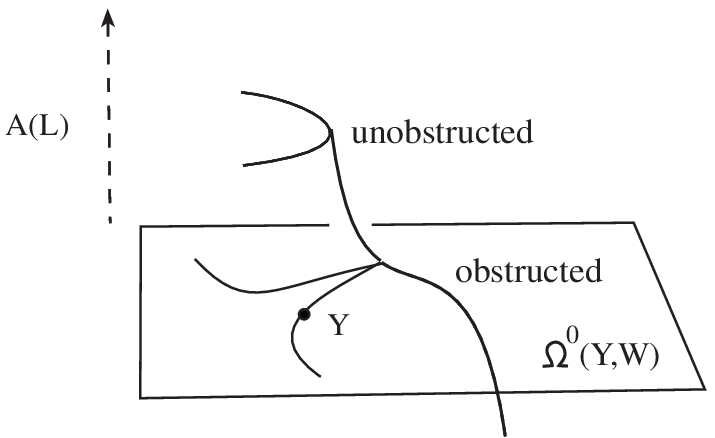}   \caption{}    \end{center}
   \end{figure}

In each slice $A$,  we are deforming along normal vector fields by the connection $\A$ (which is a perturbation of the background connection $\A_{0}$).    To get  compactness we can  cut down  the moduli space by an additional equation as in Seiberg-Witten theory, i.e.
\begin{equation}
\Phi : \Omega^{0}(Y,W) \times \CA(L)\to \Omega^{0}(Y,W) \times \Omega^{2} (Y,i\R)
\end{equation}
\begin{equation}
 \begin{array}{c}
 \Di_{\A}(v)=0 \\
F_{A}=  \sigma(v,v)
 \end{array}
  \end{equation}
  
\noindent where $F_{A}$ is the curvature of the connection $A=A_0 +a$ in $L$. For convenience we will replace the second equation above by $*F_{A}=  \mu(v,v)$,  where $*$ is the star operator on $Y$ and $\mu(v,v)=*\sigma(v,v)$.  Now we proceed exactly as in the Seiberg-Witten theory of $3$-manifolds (e.g. \cite{c}, \cite{ma}, \cite{w}). To obtain smoothness on the zeros  of $\Phi $, we  perturb the equations by $1$-forms $\delta \in \Omega^{1}(Y)$. This gives a new  parametrized equation $\Phi =0$ 
\begin{equation}
\Phi : \Omega^{0}(Y,W) \times \CA(L) \times  \Omega^{1} (Y) \to \Omega^{0}(Y,W)\times \Omega^{1} (Y,i\R)
\end{equation}
\begin{equation}
  \begin{array}{ccc}
 \Di_{A}(v) &=& 0\\
 *F_{A} + i \delta &= & \mu (v,v)
\end{array}
 \end{equation}

\noindent It follows that  in particular $\delta$ must be coclosed. $ \Phi $ has a  linearization:
\begin{equation}
D\Phi_{(v_0,A_{0},0)}: \Omega^{0}(Y,W) \times \Omega^{1}(Y,i\R) \times  \Omega^{1} (Y) 
\to \Omega^{0}(Y,W)\times \Omega^{1} (Y,i\R)
\end{equation}
$$D\Phi_{(v_0,A_{0},0)}\;(v,a,\delta)=(\Di_{A_{0}}(v) + a. v_{0}, \; *da + i\delta -2 \mu(v_{0},v))$$

\noindent We see that  $\Phi^{-1}(0)$ is smooth and the projection $\Phi^{-1}(0) \to \Omega^1 (Y)$ is onto, so by Sard's theorem for  a generic choice of $\delta$ we can make $\Phi_{\delta}^{-1}(0)$  smooth, where  $\Phi_{\delta}(v,A)=\Phi (v,A,\delta)$. The normal bundle of $Y$ has a complex structure, so  the  gauge group  $\CG(Y)=Map (Y, S^1)$ acts on the solution set  $\Phi^{-1}_{\delta}(0) $, and makes the quotient $\Phi^{-1}_{\delta }(0) / \CG(Y)$ a  smooth zero-dimensional manifold. This is because  the infinitesimal   action of $\CG(Y)$ on the complex 
$\Phi_{\delta}: \Omega^{0}(Y,W)  \times  \CA(L)  \to \Omega^{0}(Y,W)\times \Omega^{1} (Y,i\R)$
 is given by the  map 
$$\Omega^{0}(Y,i\R)\stackrel{G}{\longrightarrow} \Omega^{0}(Y,W)\times \Omega^{1}(Y,i\R)  $$
where $G( f )= ( fv_{0}, df)$.
So after dividing by $\CG$,  tangentially the complex $\Phi_{\delta}$ becomes
\begin{equation}
\Omega^{0}(Y; i\R)\stackrel{G}{\longrightarrow}  \Omega^{0}(Y,W)\times \Omega^{1}(Y,i\R)  
\to \Omega^{0}(Y,W)\times \Omega^{1} (Y,i\R)/G
\end{equation}
 Hence  the index of this complex is the sum of the indicies of the Dirac operator 
 $\Di_{A_{0}}: \Omega^{0}(Y,W)\to \Omega^{0}(Y,W)$ (which is zero), and the index of the following complex 
 \begin{equation}
 \Omega^{0}(Y,i\R)\times  \Omega^{1}(Y,i\R)  
\to \Omega^{0}(Y,i\R)\times \Omega^{1} (Y;i\R)
\end{equation}
given by $ (f,a)\mapsto (d^{*}(a), df+ *da) $, which is also zero,  since $Y^3$ has zero Euler characteristic.
  Furthermore $\Phi^{-1}_{\delta}(0) / \CG(Y)$ is compact and oriented (by Seiberg-Witten theory). Hence we can assign a number $\lambda_{\varphi}(Y,j)$ to a complex  associative submanifold $Y$ in a $G_2$ manifold . Here we don't worry about metric dependence of $\lambda_{\varphi}(Y,j)$,  since we have a fixed background metric induced  from the $G_2$ structure. In particular such a $Y$  moves in an unobstructed way after perturbation, not in the  sections  of its normal bundle $\Omega^{0}(Y,W)$,  but in  $\Omega^{0}(Y,W)\times \CA(L)$ that lies over it.   So we have:


 {\Thm Let $(Y,j)$   be a complex associative submanifold of a $G_2$ manifold $ (M,\varphi )$. Deformation space of $Y$, given by equations (32),  can be perturbed to a zero dimensional compact smooth oriented manifold, so we can associate a number $\lambda_{\varphi} (Y,j)\in \Z$. }

 \section{\bf $Spin^{c}(4)$ Structures (general case) }

We can generalize the deformations in Section 8 without  the ``complex associative" hypothesis (i.e. without the assumption that $E$  reduces to  a line bundle). We can do this is by  coupling the Dirac equation (25) with a curvature equation similar to the second equation of  (32)   we get a  generalized  version of Seiberg-Witten  theory \cite{f} (so called $SO(3)$-monopoles),  in which case we get a similar deformation picture as in Section 8 (Figure 2),  but since technically this theory is much more complicated,  instead we will put  a slight topological restriction on $M$, which will bring $Spin^{c}(4)$ structures into the picture and allow us to do the deformations in the associated bundles in an easier way by using only the standard Seiberg-Witten theory.  Also, we will first work globally  on $\tilde{M}_{\varphi }$ rather than studying deformations locally at each  submanifold.   Then at each  given associative submanifold $Y\subset M$, we can localize this  global deformations to get local deformations.

 \vspace{.1in} 

 Let $(M,\varphi)$ be a $G_2$  manifold. Now suppose the  $SO(4)$ bundle $\CP(\V) \to \tilde{M}$ lifts to a $Spin^{c}(4)$ bundle
 ${\bf V}$, call such a $(M,\varphi)$ a {\it $G_2$ manifold with $Spin^{c}(4)$ structure}. Since 
 $$Spin^{c}(4)= (SU(2)\times SU(2) \times S^{1})/{\Z}_{2}$$
a lifting exists if the Stiefel-Whitney class $w_{2}(\V)\in H^{2}(\tilde{M},{\Z}_{2})$ is  an integral class, and  the number of different liftings,  i.e. $Spin^{c}$ structures,  on $\V$ is given by $H^{2}(\tilde{M},\Z)$. These facts follow from the fibrations:
  $$S^1\to Spin^{c}(4) \to SO(4)\to K(\Z,2)\to BSpin^{c}(4)\to BSO(4)\to K(\Z,3)$$
 
$\bf {V}$ gives the  following  vector bundles by the representations,
for  $ [  q, \lambda,t ] \in Spin^{c}(4)$: 
\begin{equation}
\begin{array}{lcc}
\V \; :\;\; \;\;\;\; &y  \mapsto q y\lambda^{-1}   & \hspace{1in}  \\
 \V^{+}\; : \;\; &y  \mapsto qy t^{-1} &   \\
 \V^{-}\; : \;\; &y \mapsto  \lambda  y  t ^{-1} &\\
 ad  (\V^{+}) \; : \;\; &y \mapsto  q  y  q^{-1} &\\
  ad   (\V^{-})\; : \;\;\; &y \mapsto  \lambda  y  \lambda^{-1} &\\
  \L \;: \;\;\; &y \mapsto  y  t^{2} &
 \end{array}
 \end{equation}
   $\V^{\pm}$ are pair of $\C^{2}$ bundles that reduce to  $\R^{3}$ bundles $ad(\V^{\pm})$, and $\L$ is a complex line bundle, and  $   \V_{\C} =  \V^{+}\otimes_{\C} \bar{\V}^{-}$. Clearly $\lambda_{+}(\V)= ad (\V^{+})$, and
 $\lambda_{-}(\V)= ad (\V^{-})$.   Recall that if  we restrict  to $\tilde{M}_{\varphi }$ we get an additional  identifications  $\Xi= ad(\V^{+})$.  Also the $Spin^{c}$ structure on $\V$ gives  a $Spin^{c}$ structure on the $\R^3$ bundle $ad(\V^{+})$, which can be identifed as an equivalence class of a non-vanishing vector field $\xi_{0}$  on $ad(\V^{+})$, call this the  basic section of $ad(\V^{+})$.  
 
 \vspace{.1in} 
 
 Recall that $T(\tilde{M}_{\varphi})= T_{vert}(\tilde{M}_{\varphi} ) \oplus \Xi\oplus \V$, and  $ T_{vert}(\tilde{M}_{\varphi} )\subset \Xi^{*}\otimes \V$.  Let $\xi_{0}\in \Xi$ be the basic section. We can define an  action  $\rho: T(\tilde{M}_{\varphi})\to Hom (\V^{+}\oplus \V^{-})$, i.e. 
 $$ T(\tilde{M}_{\varphi})\otimes \V^{\pm}\to \V^{\mp}$$ 
 For  $ w=a\otimes v+x+y \in T(\tilde{M}_{\varphi}) $,  $ a\otimes v \in \Xi^{*}\otimes \V $, $ x\in \Xi $, $ y\in\V $,  let $v_{0}=a(\xi_{0})v$ and 
 $$w.(z_1,z_2)=(\bar{a}vz_{2}+ \bar{x}v_{0}z_{2}+ yz_{2} , -\bar{v}az_{1}- \bar{v_{o}}xz_{1} -\bar{y}z_1)$$
 where  $ (z_{1},z_{2})\in \V^{+}\oplus \V^{-} $.  From  (38) it is easy to check that this is a Dirac action, i.e.  $\rho(w)\circ \rho(w) =-|w|^2 I$. As usual we can dualize and extend this to  an action $$\Lambda^{2}(T^{*}\tilde{M}_{\varphi})\otimes \V^{+}\to \V^{+}$$  

 \vspace{.05in} 

\noindent Also as in (22) we have the quadratic bundle map $\sigma: \V^{+}\otimes \V^{+}\to \Lambda^{2}(\Xi^{*})_{\C}$ defined  by
 \begin{equation}
 \sigma(x,y)= -\frac{1}{2}(x i \bar{y})i
\end{equation}

\vspace{.05in}

 Now let $\CA(\L)$ be the space of connections on $\L$. Let $\A_{0}$ be a given background connection on $\V$ (e.g. the one induced from the $G_2$ metric on $M$). Than any $A\in \CA(\L)$ along with $\A_{0}$ determines a connection  on ${\bf V}$ and hence connections all the above associated  bundles. Also if  we fix $A_{0}\in \CA(\L)$,  any other $A\in \CA(\L)$ can be written as $A=A_{0}+a$ where $a\in \Omega^{1}(\tilde{M}_{\varphi})$.
Hence for each  $A\in \CA(\L)$ we get the corresponding Dirac  operator $\Di_{A}(v)=\Di_{A_{0}}(v) +a.v $ on bundles $\V^{\pm}\to \tilde{M}_{\varphi }$  which is the composition: 
$$ \Omega^{0}(\tilde{M}_{\varphi },\V^{+})\stackrel{\nabla_{A}}{\longrightarrow } \Omega^{0}(\tilde{M}_{\varphi },T^{*} \tilde{M}_{\varphi } \otimes \V^{+})
\stackrel{c}{\longrightarrow }  \Omega^{0}(\tilde{M}_{\varphi },\V^{-})  $$

 We can now write generalized Seiberg-Witten equations on $\tilde{M}_{\varphi}$ in the usual way
 $$\phi: \Omega^{0}(\tilde{M}_{\varphi}, \V^{+})\times  \CA(\L)  \to   \Omega^{0}(\tilde{M}_{\varphi},  \V^{-})  \times  \Omega^{2} (\tilde{M}_{\varphi} ) $$
 \begin{equation}
 \begin{array}{c}
 \Di_{A}(v)=0 \\
F_{A} = \sigma(v,v)
 \end{array}
  \end{equation}
and proceed as in Section 8. Equations (40) imply  $F_{A}-\sigma(v,v)$ is closed. By perturbing the equations $\phi=0$ by  closed $2$-forms $ \CZ^{2}( \tilde{M}_{\varphi})$  of $ \tilde{M}_{\varphi} $,   i.e. by changing $\phi$ to     
$$\Phi : \Omega^{0}(\tilde{M}_{\varphi}, \V^{+})\times  \CA(\L) \times   \CZ^{2}( \tilde{M}_{\varphi}) \to \Omega^{0}(\tilde{M}_{\varphi},  \V^{-}) \times  \Omega^{2} (\tilde{M}_{\varphi} )  $$

\noindent where $\Phi (A,v,\delta) =   (F_{A}  - \sigma (v,v) +i \delta, \Di_{A_{0}}(v) + a.v)$ we can make  $ \Phi ^{-1}(0)$  smooth. That is,  by applying Sard's theorem to the projection $ \Phi ^{-1}(0) \to \CZ^{2}( \tilde{M}_{\varphi}) $  for  generic choice of $\delta $ we can make $\Phi _{\delta}^{-1}(0)$  smooth, where $ \Phi_{\delta}(v,A)= \Phi (v,A,\delta)$. The normal bundle of $Y$ has a complex structure, so  the  gauge group  $\CG(\tilde{M}_{\varphi})=Map (\tilde{M}_{\varphi}, S^1)$ acts on the solution set  $\Phi^{-1}(0) $, and makes quotient $\Phi^{-1}(0) / \CG(\tilde{M}_{\varphi})$ smooth. 
The linearization of  $\Phi $ (modulo the gauge group) is given by the complex 
$$\Omega^{0}(\tilde{M}_{\varphi})\to \Omega^{1}(\tilde{M}_{\varphi}) \times \Omega^{0}(\tilde{M}_{\varphi}\V^{+})  
\to  \Omega^{2} (\tilde{M}_{\varphi} ) \times \Omega^{0}(\tilde{M}_{\varphi},\V^{-})$$
which is the sum of the  Dirac operator $\Di_{A}$ with zero index,  and a complex which is not elliptic: 
$$\Omega^{0}(\tilde{M}_{\varphi})\stackrel{d}\longrightarrow  \Omega^{1}(\tilde{M}_{\varphi}) 
\stackrel{d}\longrightarrow   \Omega^{2} (\tilde{M}_{\varphi} ) $$

\noindent where the first map is the derivative of the gauge group action.    So we have:

\vspace{.1in}

{\Thm  Normal bundle of  any associative $f: Y^{3}\hookrightarrow  M$  in a $G_2$ manifold  with a $Spin^{c}(4)$ structure $(M,\varphi, c)$, pulls back the $Spin^{c}(4)$ structure $c$  by the canonical lifting  $\tilde{f} : Y\to \tilde{M}$.   The map $\tilde{f}$  also  transforms all the bundles and their Clifford multiplication information to $Y$. So  we get a pair of $\C^2$ bundles  $V^{\pm }=\tilde{f}^{*}(\V^{\pm })$  and a line bundle $L\to Y$ describing a  $Spin^{c}(4)$ structure on $\nu(Y)=\tilde{f}^{*}(\V)$ with $\nu(Y)_{\C}=V^{+}\otimes V^{-}$. Also,  $T^{*}(Y)= \tilde{f}^{*} ad(\V^{+}) $, and any  $A\in \CA(L)$ with the help of the  background connection $A_{0}$ gives  connections on $\V^{\pm}$, and so we can transform equations (40) to $Y$ and as  in Section 8  assign  an invariant $\lambda_{\varphi}(Y,c)$ to $Y$}.

  {\Rm We can also write global version of the equations given by  (37) 

\begin{equation}
 \begin{array}{c}
 d^{*}(a)=0 \\
da +d^{*} (f\tilde{\varphi})=0
 \end{array}
  \end{equation}
 }
where $\tilde{\varphi}\in \Omega^{3}(\tilde{M})$ is the pullback of the $3$-form $\varphi $ by the projection $\tilde{M}_{\varphi} \to M$.

\end{document}